\def\Datum{January 2006}
\magnification=\magstephalf
\hsize=16.2truecm
\vsize23.2truecm
\parskip4pt plus 1pt
\frenchspacing
\lineskiplimit=-2pt
\vglue6truemm

\newread\aux\immediate\openin\aux=\jobname.aux
\ifeof\aux\else\input\jobname.aux\fi\closein\aux

\parindent25pt
\newif\ifdraft\newif\ifneu
\newdimen\SkIp\SkIp=\parskip
\pageno=1

\font\sevenrm=txr scaled 700
\font\fiverm=txr scaled 500
\font\tenrsfs=rsfs10
\font\sevenrsfs=rsfs10 scaled 700
\font\fiversfs=rsfs10 scaled 500
\font\tenbf=txb


\newfam\rsfsfam
\textfont\rsfsfam=\tenrsfs
\scriptfont\rsfsfam=\sevenrsfs
\scriptscriptfont\rsfsfam=\fiversfs
\def\rsfs{\fam\rsfsfam}

 \font\eightrm=txr scaled 800

 \font\ninerm=txr scaled 900
 \font\ninesl=txsl scaled 900
 \font\ninebf=txb scaled 900
 
 \font\gross=txb scaled\magstep1
 \font\Gross=txb scaled\magstep2

\def\nline{\hfill\break}
\long\def\fussnote#1#2{{\baselineskip=9pt
     \setbox\strutbox=\hbox{\vrule height 7 pt depth 2pt width 0pt}%
     \eightrm
     \footnote{#1}{#2}}}

\def\footnoterule{\kern-3pt
         \hrule width 2 true cm
         \kern 2.6pt}

\def\Times{\li2\times\li2}

\def\mapright#1{\smash{\mathop{\hbox to 35pt{\rightarrowfill}}\limits^{#1}}}

\def\cd{{\cdot}}
\font\trm=cmr17

\font\nineti=cmmi9

\font\ninett=cmtt9

\font\ninesy=cmsy9

\def\mod{\mathop{\rm mod}\nolimits}
\def\rank{\mathop{\rm rank}\nolimits}

\def\Alpha{\hbox{\mib\char11}}

\font\tenmsb=msbm10\font\sevenmsb=msbm7\font\fivemsb=msbm5
\newfam\msbfam
\textfont\msbfam=\tenmsb\scriptfont\msbfam=\sevenmsb
\scriptscriptfont\msbfam=\fivemsb

\font\tenmib=cmmib10 \font\eightmib=cmmib8 \font\fivemib=cmmib5
\newfam\mibfam
\textfont\mibfam=\tenmib \scriptfont\mibfam=\eightmib
\scriptscriptfont\mibfam=\fivemib
\def\mib{\fam\mibfam\tenmib}

\font\tenmi=cmmi10 \font\eightmi=cmmi8 \font\fivemi=cmmi5
\newfam\mifam
\textfont\mifam=\tenmi \scriptfont\mifam=\eightmi
\scriptscriptfont\mifam=\fivemi

\font\tenfr=eufm10 \font\eightfr=eufm8 \font\fivefr=eufm5
\newfam\frfam
\textfont\frfam=\tenfr \scriptfont\frfam=\eightfr
\scriptscriptfont\frfam=\fivefr
\def\frak{\fam\frfam\tenfr}

\font\tenbfr=eufb10 \font\eightbfr=eufb8 \font\fivebfr=eufb5
\newfam\bfrfam
\textfont\bfrfam=\tenbfr \scriptfont\bfrfam=\eightbfr
\scriptscriptfont\bfrfam=\fivebfr

\def\Quot#1#2{\raise 2pt\hbox{$#1\mskip-1.2\thinmuskip$}\big/%
     \lower2pt\hbox{$\mskip-0.8\thinmuskip#2$}}
\def\dd#1{\raise1.5pt\hbox{$\re1\partial\!$}/\raise-2.5pt\hbox{$\!\partial#1$}}
\def\ddz#1{\raise1.5pt\hbox{$\re1\partial\!$}/\raise-2.5pt\hbox{$\!\partial z_{#1}$}}
\def\ddw#1{\raise1.5pt\hbox{$\re1\partial\!$}/\raise-2.5pt\hbox{$\!\partial w_{#1}$}}
\def\hol{\7{hol}}\def\aut{\7{aut}}

\def\Klein{\ninerm\textfont1=\nineti\textfont0=\ninerm\def\sl{\ninesl}
\textfont2=\ninesy\baselineskip11.2pt\def\bf{\ninebf}}

\def\Kl#1{\raise1pt\hbox{$\scriptstyle($}#1\raise1pt\hbox{$\scriptstyle)$}}

\newcount\ite\ite=1
\def\0{\global\ite=1\1}
\def\1{\item{\rm(\romannumeral\the\ite)}\global\advance\ite1}
\def\3#1{{\mib#1}}
\def\5#1{{\cal#1}}
\def\6#1{{\rsfs#1}}
\def\7#1{\mathop{\frak#1\re2}\nolimits}
\def\8#1#2#3#4{\big(\raise-2pt\hbox{\rlap{\raise6pt%
  \hbox{$\scriptstyle#1\re1#2$}}\hbox{$\scriptstyle#3\re2#4$}}\big)}\
\def\9#1{\rlap{\hskip.2pt$#1$}\rlap{\hskip.4pt$#1$}#1}

\def\re#1{\mskip#1mu}
\def\li#1{\mskip-#1mu}

\def\tilde{\widetilde}

\def\<{\;\;\Longleftrightarrow\;\;}
\def\To#1#2{(\romannumeral#1) $\Longrightarrow$ (\romannumeral#2)}

\def\epsilon{\varepsilon}
\def\hat{\rlap{\raise-1.5pt\hbox{\hskip2.5pt\trm\char94}}}

\def\phi{\varphi}
\def\steil#1{\hbox{\rm~~#1~~}}
\def\Steil#1{\hbox{\rm\quad #1\quad}}
\def\p{p.\nobreak\hskip2pt}

\def\Re{\hbox{\rm Re}\re1}

\def\MOP#1{\expandafter\edef\csname #1\endcsname{%
   \mathop{\hbox{\rm #1}}\nolimits
    }}

\MOP{Aut}\MOP{ad}\MOP{GL}\MOP{SU}\MOP{SO}\MOP{tr}\MOP{Sp}
\MOP{SL}\MOP{id}\MOP{Ad}\MOP{PSL}\MOP{U}\MOP{O}\MOP{End}

\def\sqr#1#2{{\,\vcenter{\vbox{\hrule height.#2pt\hbox{\vrule width.#2pt
height#1pt \kern#1pt\vrule width.#2pt}\hrule height.#2pt}}\,}}
\def\qed{\hfill\ifmmode\sqr66\else$\sqr66$\par\fi\rm}

\def\One{{1\kern-3.8pt 1}}
\def\one{{1\kern-3.1pt 1}}

\font\tenmia=txmia
\font\sevenmia=txmia scaled 700
\font\fivemia=txmia scaled 500

\def\Bbord#1{\mathchoice
{\hbox{\tenmia\char#1}}
{\hbox{\tenmia\char#1}}
{\hbox{\sevenmia\char#1}}
{\hbox{\fivemia\char#1}}
}

\def\CC{\Bbord{131}}

\def\NN{\Bbord{142}}

\def\RR{\Bbord{146}}
\def\ZZ{\Bbord{154}}

\newcount\nummer
\newcount\parno\parno=0

\def\KAP#1#2{\par\bigbreak\global\advance\parno by1%
\def\test{#1}\ifx\test\empty\else%
\expandafter\let\csname#1\endcsname=\relax%
\immediate\write\aux{\def\csname#1\endcsname{\the\parno}}%
\expandafter\xdef\csname#1\endcsname{\the\parno}\fi%
{\vskip3pt\noindent\gross\the\parno. #2\hfil}\vskip3pt\rm\nummer=0\nobreak}

\newwrite\aux

\def\Randmark#1{\vadjust{\vbox to 0pt{\vss\hbox to\hsize%
{\fiverm\hskip\hsize\hskip1em\raise 2.5pt\hbox{#1}\hss}}}}

\def\PrN{\the\parno.\the\nummer}

\def\Write#1#2{\global\advance\nummer1\def\test{#1}\ifx\test\empty\else%
\ifdraft\Randmark{#1}\fi\expandafter\let\csname#1\endcsname=\relax%
\immediate\write\aux{\def\csname#1\endcsname{\PrN}}%
\expandafter\xdef\csname#1\endcsname{\PrN}\fi#2}

\def\Num#1{\Write{#1}{\PrN}}
\def\Leqno#1{\Write{#1}{\leqno(\PrN)}}

\def\Proposition#1{{\smallbreak\noindent\bf\Num{#1} Proposition.~}\parskip0pt\sl}
\def\Lemma#1{{\par\noindent\bf\Num{#1} Lemma.~}\parskip0pt\sl}
\def\Corollary#1{{\smallbreak\noindent\bf\Num{#1} Corollary.~}\parskip0pt\sl}

\def\Remark#1{{\smallbreak\noindent\bf\Num{#1} Remark.~}\parskip0pt\rm}
\def\Construction#1{{\smallbreak\noindent\bf\Num{#1} Construction.~}\parskip0pt\rm}
\def\Example#1{{\smallbreak\noindent\bf\Num{#1} Example.~}}
\def\Definition#1{{\smallbreak\noindent\bf\Num{#1} Definition.~}}
\long\def\Proof{\smallskip\noindent\sl Proof.\parskip\SkIp\rm~~}
\def\Formend{\par\parskip\SkIp\rm}

\def\ruf#1{{{\expandafter\ifx\csname#1\endcsname\relax\xdef\flAG{}%
\message{*** #1 nicht definiert!! ***}\ifdraft\Randmark{#1??}\fi\else%
\xdef\flAG{1}\fi}\ifx\flAG\empty{\bf??}\else\rm\csname#1\endcsname\fi}}
\def\Ruf#1{{\rm(\ruf{#1})}}

\newcount\lit\lit=1
\def\Ref#1{\item{\the\lit.}\expandafter\ifx\csname#1ZZZ\endcsname\relax%
\fi%
\expandafter\let\csname#1\endcsname=\relax%
\immediate\write\aux{\def\csname#1\endcsname{\the\lit}}\advance\lit1\ifdraft\Randmark{#1}\fi}

\def\Lit#1{\expandafter\gdef\csname#1ZZZ\endcsname{1}[\ruf{#1}]}
\def\LIT#1#2{\expandafter\gdef\csname#1ZZZ\endcsname{1}[\ruf{#1}, #2]}

\rm

\immediate\openout\aux=\jobname.aux

\ifdraft
\footline={\hss\sevenrm \Datum\hss}
\else\nopagenumbers\fi

\headline={\ifnum\pageno>1\sevenrm\ifodd\pageno Homogeneous Levi 
degenerate CR-manifolds
\hss{\tenbf\folio}\else{\tenbf\folio}\hss{Fels-Kaup}
\fi\else\hss\fi}

\centerline{\Gross Homogeneous Levi degenerate}
\bigskip
\centerline{\Gross CR-manifolds in dimension 5}
\bigskip\bigskip
\centerline{Gregor Fels\quad and\quad Wilhelm Kaup} 
\bigskip
{\parindent0pt\footnote{}{\ninerm 2000 Mathematics 
Subject Classification: 32M17, 32M25, 32V25.
}}

\KAP{Introduction}{Introduction} Levi nondegenerate real-analytic
hypersurfaces $M\subset\CC^{n}$ are well understood due to the seminal
results of Tanaka \Lit{TANA} and Chern-Moser\Lit{CHMO}, where a
complete set of {\sl local} invariants for $M$ has been defined. At
the other extreme are those hypersurfaces which are locally
CR-equivalent to a product $M'\times\CC$, (take the real hyperplanes
in $\CC^n$ as a simple example).  In-between there is a huge and much
less understood class of CR-hypersurfaces which are neither Levi
nondegenerate nor locally direct products as above.  A convenient
description of that class of CR-manifolds can be given in terms of
$k$-nondegeneracy in the sense of \Lit{BAHR} with $k\ge 2,$ see
\Lit{BERO} for further details.  The (uniformly) 2-nondegenerate
real-analytic hypersurfaces are the main focus of this paper.

Clearly, $5$ is the lowest dimension of $M$ (that is, $n=3$) for which
Levi degenerate hypersurfaces can occur that are not locally
equivalent to $M'\times\CC$ for any CR-manifold $M'$. Of a particular
interest are those manifolds that are Levi degenerate at every
point. A well studied example of this type is the tube over the future
light cone in 3-dimensional space-time, more precisely
$$\5M:=\{z\in\CC^{3}: (\Re z_{1})^{2}+(\Re z_{2})^{2}=(\Re
z_{3})^{2},\;\Re(z_{3})>0\}\,,$$ compare e.g. \Lit{EBFT}, \Lit{FEKA},
\Lit{KAZT} and \Lit{SEVL}. $\,\5M$ is nondegenerate in a higher order
sense (2-nondegenerate, to be precise) and is also homogeneous (a
$7$-dimensional group of complex affine transformations acts
transitively on it). All further examples of locally homogeneous
$2$-nondegenerate CR-manifolds of dimension 5 known so far finally
turned out to be locally CR-isomorphic to the light cone tube $\5M$,
compare e.g. \Lit{EBEN}, \Lit{KAZT}, \Lit{GAME}, \Lit{FEKA}, and the
belief arose that there are no others. The main objective of this
paper is the presentation of an infinite family of
homogeneous $2$-nondegenerate CR-manifolds of dimension 5 which are
pairwise locally CR-inequivalent. In the forthcoming paper \Lit{FEKP}
it will be shown by purely Lie algebraic methods that our list of
examples is complete in the following sense: {\sl Every locally
homogeneous $2$-nondegenerate real-analytic CR-manifold of dimension 5
is locally CR-equivalent to one of our examples.}

A main point in our discussion is to show that the provided examples
(Examples \ruf{EI} -- \ruf{EV}) are mutually CR-nonequivalent. It
should be noted that even in case of Levi nondegenerate hypersurfaces
it is in general quite hard to decide whether these are locally
CR-equivalent at selected points, although this is possible in
principle due to the existence of a complete system of invariants
\Lit{CHMO}. In the Levi degenerate case much less is known and
suitable invariants have to be created ad hoc. In our particular case
of locally homogeneous CR-manifolds $M$ we use for every $a\in M$ as
invariant the Lie algebra isomorphy type of $\hol(M,a)$, which by
definition is the real Lie algebra of all germs at $a$ of local
infinitesimal CR-transformations defined in suitable neighbourhoods of
$a\in M$, see Section \ruf{Preliminaries} for details. It is meanwhile
well known that for the tube $\5M$ over the future light cone in
$\RR^{3}$ the Lie algebras $\hol(\5M,a)$ have dimension 10 and are
isomorphic to the simple Lie algebra $\7{so}(2,3)$. It turns out that
for all our examples $M\ne\5M$ of $2$-nondegenerate locally
homogeneous CR-manifolds the Lie algebras $\hol(M,a)$ are solvable of
dimension 5 and are mutually non-isomorphic for different
$M$. Clearly, this procedure needs a feasible way to compute the Lie
algebras $\hol(M,a)$ explicitly. One of our major goals is to provide
such a way at least for manifolds of tube type.

The paper is organized as follows. After recalling some necessary
preliminaries in Section \ruf{Preliminaries} we discuss in Section
\ruf{Tube} tube manifolds $M=F\oplus i\RR^{n}\subset\CC^{n}$ over
real-analytic submanifolds $F\subset\RR^{n}$. It turns out that the
CR-structure of $M$ is closely related to the real-affine structure of
the base $F$. For instance, the Levi form of $M$ is essentially the
sesquilinear extension of the second fundamental form of the
submanifold $F\subset\RR^{n}$.  Generalizing the notion of the second
fundamental form we define higher order invariants for $F$ (see
\ruf{HN}).  In the uniform case these invariants precisely detect the
$k$-nondegeneracy of the corresponding CR-manifold $M=F\oplus i\RR^n$.
It is known that the (uniform) $k$-nondegeneracy of a real-analytic
CR-manifold $M$ together with minimality ensures that the Lie algebras
$\hol(M,a)$ are finite-dimensional.  For submanifolds
$F\subset\RR^{n}$, homogeneous under a group of affine
transformations, a simple criterion for the $2$-nondegeneracy of the
associated tube manifold $M$ is given in Proposition \ruf{US}. We
close this section with some general remarks concerning polynomial
vector fields.

In Section \ruf{cones} these results are applied to the case where $F$
is conical in $\RR^n$, that is, locally invariant under dilations
$z\mapsto tz$ for $t$ near $1\in\RR$.  In this case $M$ is always {\sl
Levi degenerate.}  Assuming that $\hol(M,a)$ is finite dimensional
(equivalently, $M$ is holomorphically non-degenerate and minimal) we
develop some basic techniques which enable us to compute $\hol(M,a)$.
Our first observation is that under the finiteness assumption
$\hol(M,a)$ consists only of polynomial vector fields and carries a
natural graded structure, see Lemma \ruf{DC}.  We prove further (under
the same assumptions, see Proposition \ruf{ZB}) that local
CR-equivalences between two such tube manifolds are always rational
(even if these manifolds are not real-algebraic).  If all homogeneous
parts in $\hol(M,a)$ of degree $\ge2$ vanish, this result can be
further strengthened, see Proposition \ruf{LO}.  We close the section
with an investigation of the tubes $M^\alpha_{p,q}$ over the cones
$F^\alpha_{p,q}:=\{x\in \RR^{p{+}q}_+: \sum \epsilon_j x_j^\alpha=0\}$
with $p$ positive and $q$ negative $\epsilon_j$'s and $\alpha\in
\RR^{*}$. We show how the preceeding results can be used to determine
explicitly all $\hol(M^\alpha_{p,q},a)$ for $\alpha\ge2$ an arbitrary
integer.

In Section \ruf{Some}, we study homogeneous CR-submanifolds
$M^{\phi}\subset\CC^{n}$ depending on the choice of an endomorphism
$\phi\in \End(\RR^n)$ together with a cyclic vector $a\in\RR^{n}$ for
$\phi$ in the following way: The powers $\phi^0,\phi^1,\ldots,\phi^d$
($2\le d+1<n$) span an abelian Lie algebra $\7h^\phi$ and, in turn,
give the cone $F^{\phi}=\exp(\7h^{\phi})a\subset \RR^n$.  For
`generic' $\phi$ the corresponding tube manifolds
$M^\phi=F^{\phi}\oplus\,i\RR^{n}$ are 2-nondegenerate and of
CR-dimension $d+1$.  In this case the local invariants
$\hol(M^\phi,a)$ are explicitly determined (Proposition \ruf{UH}).
Further, again for $\phi$ in general position the tube manifold
$M^\phi$ is simply connected and has trivial stability group at every
point. As a consequence, the manifolds $M=M^{\phi}$ of this type have
the following remarkable property: {\sl Every homogeneous
(real-analytic) CR-manifold locally CR-equivalent to $M$ is globally
CR-equivalent to $M$} (see \ruf{UQ}). Specialized to $n=3$ we get in
Section \ruf{Examples} an infinite family of homogeneous
$2$-nondegenerate hypersurfaces in $\CC^{3}$. These are all tubes over
linearly homogeneous cones in $\RR^{3}$. Besides these, there is just
one other example, namely the tube over the twisted cubic in $\RR^{3}$
(which is not a cone but is affinely homogeneous).

\KAP{Preliminaries}{Preliminaries} In the following let $E$ always be
a complex vector space of finite dimension and $M\subset E$ a locally
closed connected real-analytic submanifold (unless stated otherwise).
Due to the canonical identifications $T_aE=E,$ for every $a\in M$ we
consider the tangent space $T_{a}M$ as an $\RR$-linear subspace of
$E$. Set $H_{a}M:=T_{a}M\cap iT_{a}M$. The manifold $M$ is called a
{\sl CR-submanifold} if the dimension of $H_{a}M$ does not depend on
$a\in M$. This complex dimension is called the {\sl CR-dimension} of
$M$ and $H_{a}M$ is called the {\sl holomorphic tangent space} at $a$,
compare \Lit{BERO} as general reference for CR-manifolds. Given a
further real-analytic CR-submanifold $M'$ of a complex vector space
$E'$ a smooth mapping $ h:M\to M'$ is called {\sl CR} if for all $a\in
M$ the differentials $d h_{a}:T_{a}M\to T_{ h a}M'$ map the
corresponding holomorphic tangent spaces in a complex linear way to
each other. A vector field on $M$ is a mapping $f:M\to E$ with
$f(a)\in T_{a}M$ for all $a\in M$. For better distinction we also
write $\xi=f(z)\dd z$ instead of $f$ and $\xi_{a}$ instead of $f(a)$,
compare the convention (2.1) in \Lit{FEKA}. The real-analytic vector
field $\xi=f(z)\dd z$ is called an {\sl infinitesimal
CR-transformation} on $M$ if the corresponding local flow consists of
CR-transformations.

For simplicity and without essential loss of generality we always
assume that the CR-submanifold $M$ is {\sl generic} in $E$, that is,
$E=T_{a}M+iT_{a}M$ for all $a\in M$. Given an infinitesimal
CR-transformation $f(z)\dd z$ on a generic CR-submanifold $M$ the map
$f$ extends uniquely to a holomorphic map $f:U_f\to E$ in an open
neighborhood $U_f$ of $M$ in $E$, (see \Lit{ANHI} or 12.4.22 in
\Lit{BERO}). Let us denote by $\hol(M)$ the space of all such vector
fields, which is a real Lie algebra with respect to the usual
bracket. Further, for every $a\in M$ we denote by $\hol(M,a)$ the Lie
algebra of all germs of infinitesimal CR-transformations defined in
arbitrary open neighbourhoods of $a\in M$. Since $M$ is generic in $E$
the space $\hol(M,a)$ naturally embeds as a real Lie subalgebra of the
complex Lie algebra $\hol(E,a)$. The CR-manifold $M$ is called {\sl
holomorphically nondegenerate} at $a$ if $\hol(M,a)$ is totally real
in $\hol(E,a)$, that is, if $\hol(M,a)\cap i\hol(M,a)=0$ in
$\hol(E,a)$. This condition together with the minimality of $M$ in $E$
are equivalent to $\dim \hol(M,a)<\infty$ (see (12.5.16) in
\Lit{BERO}). Here, the CR-submanifold $M\subset E$ is called {\sl
minimal} at $a\in M$ if $T_{a}R=T_{a}M$ for every submanifold
$R\subset M$ with $a\in R$ and $H_{a}M\subset T_{a}R$.
 
By $\aut(M,a):=\{\xi\in\hol(M,a):\xi_{a}=0\}$ we denote the {\sl
isotropy subalgebra} at $a\in M$. Clearly, $\aut(M,a)$ has finite
codimension in $\hol(M,a)$. Furthermore, we denote by $\Aut(M,a)$ the
group of all {\sl germs} of real-analytic CR-isomorphisms $ h:W\to
\tilde W$ with $ h(a)=a$, where $W,\tilde W$ are arbitrary open
neighbourhoods of $a$ in $M$. It is known that every germ in
$\Aut(M,a)$ can be represented by a holomorphic map $U\to E$, where
$U$ is an open neighbourhood of $a$ in $E$, compare e.g. 1.7.13 in
\Lit{BERO}. Finally, $\Aut(M)$ denotes the group of all global
real-analytic CR-automorphisms $h:M\to M$ and $\Aut(M)_a$ its isotropy
subgroup at $a$. There is a canonical group monomorphism
$\Aut(M)_a\hookrightarrow\Aut(M,a)$ as well as an exponential map
$\exp:\aut(M,a)\to\Aut(M,a)$ for every $a\in M$.

\medskip A basic invariant of a CR-manifold is the (vector valued)
{\sl Levi form}. Its definitions found in the literature may differ by
a constant factor. Here we choose the following definition: It is
well-known that for every point $a$ in the CR-manifold $M$ there is a
well defined alternating $\RR$-bilinear map
$$\omega_{a}:H_{a}M\times H_{a}M\;\to\; E/H_{a}M$$ satisfying
$\omega_{a}(\xi_{a},\zeta_{a})=[\xi,\zeta]_{a}\steil{mod}H_{a}M,$
where $\xi,\zeta$ are arbitrary smooth vector fields on $M$ with
$\xi_{z},\zeta_{z}\in H_{z}M$ for all $z\in M$. We define the Levi
form
$$\5L_{a}:H_{a}M\times H_{a}M\;\to\; E/H_{a}M\Leqno{LE}$$ to be twice
the sesquilinear part of $\omega_{a}$. By {\sl sesquilinear} we always
mean 'conjugate linear in the first and complex linear in the second
variable', that is,
$$\5L_{a}(v,w)=\omega_{a}(v,w)+i\omega_{a}(iv,w)\qquad \hbox{for all }
v,w\in H_{a}M\,.$$ In particular, the vectors $\5L_{a}(v,v),$ $v\in
H_aM,$ are contained in $iT_{a}M/H_{a}M$, which can be identified in a
canonical way with the normal space $E/T_{a}M$ to $M\subset E$ at $a$.

Define the {\sl Levi kernel}
$$K_{a}M:=\{v\in H_{a}M:\5L_{a}(v,w)=0\steil{for all}w\in H_aM\}\,.$$
The CR-manifold $M$ is called {\sl Levi nondegenerate} at $a$ if
$K_{a}M=0$. Generalizing that, the notion of {\sl $k$-nondegeneracy}
for $M$ at $a$ has been introduced for every integer $k\ge1$ (see
\Lit{BAHR}, \Lit{BERO}). As shown in 11.5.1 of \Lit{BERO} a
real-analytic and connected CR-manifold $M$ is holomorphically
nondegenerate at $a$ (equivalently: at every $z\in M$) if and only if
$M$ is $k$-nondegenerate at $a$ for some $k\ge1$. For $k=1$ this
notion is equivalent to $M$ being Levi nondegenerate at $a\in M$.

\bigskip In Section \ruf{Some} we also need a more general notion of
(real-analytic) CR-manifold. This is a connected real-analytic
manifold $M$ together with a subbundle $HM\subset TM$ and a bundle
endomorphism $J:HM\to HM$ satisfying the following property: For every
point of $M$ a suitable open neighbourhood can be realized as (locally
closed) CR-submanifold $U$ of some $\CC^{n}$ in such a way that
$H_{z}M$ corresponds to $H_{z}U=T_{z}U\cap iT_{z}U$ and the
restriction of $J$ to $H_{z}M$ corresponds to the multiplication with
the imaginary unit $i$ on $H_{z}U$ for every $z\in U\subset M$.

\KAP{Tube}{Tube manifolds}

Let $V$ be a real vector space of finite dimension and $E:=V\,\oplus\,
iV$ its complexification. Let furthermore $F\subset V$ be a connected
real-analytic submanifold and $M:=F\oplus iV\subset E$ the
corresponding tube manifold. $M$ is a generic CR-submanifold of $E$
invariant under all translations $z\mapsto z+iv$, $v\in V$. In case
$V'$ is another real vector space of finite dimension, $E'$ is its
complexification, $F'\subset V'$ a real-analytic submanifold and
$\phi:V\to V'$ an affine mapping with $\phi(F)\subset F'$, then
clearly $\phi$ extends in a unique way to a complex affine mapping
$E\to E'$ with $\phi(M)\subset M'$. It should be noted that higher
order real-analytic maps $\psi:F\to F'$ also extend locally to
holomorphic maps $\psi:U\to E'$, $U$ open in $E$. But in contrast to
the affine case we have in general $\psi(M\cap U)\not\subset M'$. We
may therefore ask how the CR-structure of $M$ is related to the real
affine structure of the submanifold $F\subset V$.

For every $a\in F$ let $T_{a}F\subset V$ be the {\sl tangent space}
and $N_{a}F:=V/T_{a}F$ the {\sl normal space} to $F$ at $a$. Then
$T_{a}M=T_{a}F\oplus iV$ for the corresponding tube manifold $M$, and
$N_{a}F$ can be canonically identified with the normal space
$N_{a}M=E/T_{a}M$ of $M$ in $E$. Define the map $\ell_{a}:T_{a}F\times
T_{a}F\to N_{a}F$ in the following way: For $v,w\in T_{a}F$ choose a
smooth map $f:V\to V$ with $f(a)=w$ and $f(x)\in T_{x}F$ for all $x\in
F$ (actually it suffices to choose such an $f$ only in a small
neighborhood of $a\,$). Then put
$$\ell_{a}(v,w):=f'(a)(v)\;\mod\;T_{a}F\,,\Leqno{EL}$$ where the
linear operator $f'(a)\in\End(V)$ is the derivative of $f$ at $a$. One
shows that $\ell_a$ does not depend on the choice of $f$ and is a
symmetric bilinear map. We mention that if $V$ is provided with the
flat Riemannian metric and $N_aF$ is identified with $T_aF^\perp$ then
$\ell$ is nothing but the second fundamental form of $F$ (see the
subsection II.3.3 in \Lit{SAKA}). The form $\ell_{a}$ can also be read
off from local equations for $F$, more precisely, suppose that
$U\subset V$ is an open subset, $W$ is a real vector space and $h:U\to
W$ is a real-analytic submersion with $F=h^{-1}(0)$. Then the
derivative $h'(a):V\to W$ induces a linear isomorphism $N_{a}F\cong W$
and modulo this identification $\ell_{a}$ is nothing but the second
derivative $h''(a):V\times V\to W$ restricted to $T_{a}F\times
T_{a}F$.

By
$$K_{a}F:=\{w\in T_{a}F:\ell_{a}(v,w)=0\steil{for all}v\in T_{a}F\}$$
we denote the {\sl kernel} of $\ell_{a}$. The manifold $F$ is called
{\sl nondegenerate} at $a$ if $K_{a}F=0$ holds. The following
statement follows directly from the definition of $\ell_a$:

\Lemma{CR} Suppose that $\phi\in\End(V)$ satisfies $\phi(x)\in T_{x}F$
for all $x\in F$. Then $\phi(a)\in K_{a}F$ if and only if
$\phi\big(T_{a}F\big)\subset T_{a}F$.\Formend

\noindent Lemma \ruf{CR} applies in particular for $\phi=\id$ in case
$F$ is a {\sl cone,} that is, $rF=F$ for all real $r>0$. More
generally, we call the submanifold $F\subset V$ {\sl conical} if $x\in
T_{x}F$ for all $x\in F$. Then $\RR a\subset K_{a}F$ holds for all
$a\in F$.

\medskip

In the remaining part of this section we explain how the CR-structure
of the tube manifold $M$ is related to the real objects $\ell_a,\,
TF,\, KF,\, K^k\!F$. In general it is quite hard to check whether a
given CR-manifold $M$ is $k$-nondegenerate at a point $a\in M$. For
affinely homogeneous tube manifolds, however, there are simple
criteria, see Propositions \ruf{UR} and \ruf{US}. We start with some
preparations.

For every $a\in F\subset M$ $$H_{a}M=T_{a}F\oplus iT_{a}F\,\subset\,
E$$ is the holomorphic tangent space at $a$, and $E/H_{a}M$ can be
canonically identified with $N_{a}F\oplus iN_{a}F$. It is easily seen
that the Levi form $\5L_{a}$ of $M$ at $a$, compare \Ruf{LE}, is
nothing but the sesquilinear extension of the form $\ell_{a}$ from
$T_{a}F\times T_{a}F$ to $H_{a}M\times H_{a}M$. In particular,
$$K_{a}M=K_{a}F\oplus iK_{a}F$$ is the {\sl Levi kernel} of $M$ at
$a$. In case the dimension of $K_{a}F$ does not depend on $a\in F$
these spaces form a subbundle $KF\subset TF$. In that case to every
$v\in K_{a}F$ there exists a smooth function $f:V\to V$ with $f(a)=v$
and $f(x)\in K_{x}F$ for all $x\in F$, i.e., the tangent vector $v$
extends to a smooth section in $KF$. In any case, let us define
inductively linear subspaces $K_{a}^{k}F$ of $T_{a}F$ as follows:

\Definition{HN}  For every real-analytic submanifold $F\subset V$ and
every $k\in\NN$ put \1 $K_{a}^{0}F:=T_{a}F$, \1 $K_{a}^{k+1}\!F$ is
the space of all vectors $v\in K_{a}^{k}F$ such that there is a smooth
mapping $f:V\to V$ with $f'(a)(T_{a}F)\subset K_{a}^{k}F$, $f(a)=v$
and $f(x)\in K_{x}^{k}F$ for all $x\in F$.

\medskip\noindent It is clear that $K_{a}^{1}F=K_{a}F$ holds. Let us
call $F$ {\sl uniformly degenerate} if for every $k\in\NN$ the
dimension of $K_{a}^{k}F$ does not depend on $a\in F$. In that case
for every $v\in K_{a}^{k}F$ the outcome of the condition
$f'(a)(T_{a}F)\subset K_{a}^{k}F$ in (ii) does not depend on the
choice of the smooth mapping $f:V\to V$ satisfying $f(a)=v$ and
$f(x)\in K_{x}^{k}F$ for all $x\in F$. For instance, $F$ is of uniform
degeneracy if $F$ is locally affinely homogeneous, that is, if there
exists a Lie algebra $\7a$ of affine vector fields on $V$ such that
every $\xi\in\7a$ is tangent to $F$ and such that the canonical
evaluation map $\7a\to T_{a}F$ is surjective for every $a\in
F$. Clearly, if $F$ is locally affinely homogeneous in the above sense
then the corresponding tube manifold $M=F\oplus iV$ is locally
homogeneous as CR-manifold.

\medskip We identify every smooth map $f:V\to V$ with the
corresponding smooth vector field $\xi=f(x)\dd x$ on $V$. Our
computations in the following are considerably simplified by the
obvious fact that every smooth vector field $\xi$ on V has a unique
smooth extension to $E$ that is invariant under all translations
$z\mapsto z+iv$, $v\in V$. In case $\xi$ is tangent to $F\subset V$
the extension satisfies $\xi_{z}\in H_{z}M$ for all $z\in M$. 

In case $F\subset V$ is uniformly degenerate in a neighbourhood of
$a\in F$ we call $F$ {\sl $k$-nondegenerate} at $a$ if $K_{a}^{k}F=0$
and $k$ is minimal with respect to this property. As a consequence of
\Lit{KAZT}, compare the last lines therein, we state:

\Proposition{UR} Suppose that $F$ is uniformly degenerate in a
neighbourhood of $a\in F$. Then the corresponding tube manifold
$M=F\oplus iV$ is $k$-nondegenerate as CR-manifold at $a\in M$ if and
only if $F$ is $k$-nondegenerate at $a$ in the above defined
sense.\Formend

\Corollary{AF} Suppose that $F$ is uniformly degenerate, $\dim(F)\ge2$
and $K_{x}F=\RR x$ holds for every $x\in F$. Then $F$ is
$2$-nondegenerate at every point.

\Proof The condition $K_{x}F=\RR x$ implies $0\not\in F$, otherwise
the uniform degeneracy of $F$ would be violated. The map $f=\id$ has
the property $f(x)\in K_xF$ for every $x\in F$. Hence, the relation
$f'(x)(T_{x}F)=T_{x}F\not\subset K_{x}F$ implies $x\notin K_{x}^{2}F$
and thus $K_{x}^{2}F=0$.\qed

\medskip For locally affinely homogeneous submanifolds
$F\subset V$ the spaces $K_{a}^{k}F$ can easily be computed.

\Proposition{US} Suppose that $\7a$ is a linear space of affine vector
fields on $V$ such that every $\xi\in\7a$ is tangent to $F$ and the
canonical evaluation mapping $\7a\to T_{a}F$ is a linear
isomorphism. Then, given any $k\in\NN$, the vector $v\in K_{a}^{k}F$
is in $K_{a}^{k+1}\!F$ if and only if for every $\xi=h(x)\dd x\in\7a$
the relation $\lambda^h(v)\in K_{a}^{k}F$ holds, where
$\lambda^h:=h-h(0)$ is the linear part of $\xi$.

\Proof By the implicit function theorem, there exist open
neighbourhoods $Y$ of $0\in\7g$ and $X$ of $a\in M$ such that
$g(y):=\exp(y)a$ defines a diffeomorphism $g:Y\to X$. Define the
smooth map $f:X\to V$ by $f(g(y))=\mu_{y}(v)$, where $\mu_{y}$ is the
linear part of the affine transformation $\exp(y)$. Then $f(a)=v$ and
$f(x)\in K_{x}^{k}F$ for every $x\in X$. The claim now follows from
$f'(a)\big(g'(0)\xi\big)=\lambda^h(v)$ for every $\xi\in\7g$ and
$\lambda^h$ the linear part of $h$.\qed

\medskip It is easily seen that a necessary condition for $M$ being
minimal as CR-manifold is that $F$ is not contained in an affine
hyperplane of $V$. A sufficient condition is that the image of the
form $\ell_{a}$, see \Ruf{EL}, spans the full normal space
$N_{a}F$ at every $a\in F$.

\smallskip For later use in Proposition \ruf{UQ} we state

\Lemma{ZI} Suppose that $F\subset V$ is a submanifold such that for
every $c\in V$ with $c\ne0$ there exists a linear transformation
$\lambda\in\GL(V)$ with $\lambda(F)=F$ and $\lambda(c)\ne c$ (this
condition is automatically satisfied if $F$ is a cone). Then for
$M=F\oplus iV$ the CR-automorphism group $\Aut(M)$ has trivial center.

\Proof Let an element in the center of $\Aut(M)$ be given and let
$h:U\to E$ be its holomorphic extension to an appropriate connected
open neighbourhood $U$ of $M$. Since $h$ commutes with
every translation $z\mapsto z+iv$, $v\in V$, it is a translation
itself: Indeed, for $a\in F$ fixed and $c:=h(a)-a$ the translation
$\tau(z):=z+c$ coincides with $h$ on $a+iV$ and hence on $U$ by the
identity principle. In case $c\ne0$ choose $\lambda\in\GL(V)$ as in
the above assumption. Then $h$ commutes with $\lambda\in\Aut(M)$ and
$\lambda(c)=c$ gives a contradiction, showing $h(z)\equiv z$.\qed

\medskip For fixed complex vector space $E$ as above let us denote by
$\7P$ the complex Lie algebra of all polynomial holomorphic vector
fields $f(z)\dd z$ on $E$, that is, $f:E\to E$ is a polynomial
map. Then $\7P$ has the $\ZZ$-grading
$$\7P=\bigoplus_{k\in\ZZ}\7P_{k}\,,\qquad[\7P_{k},\7P_{l}]\subset
\7P_{k+l}\,,\Leqno{HJ}$$ where $\7P_{k}$ is the $k$-eigenspace of
$\ad(\delta)$ for the Euler field $\delta:=z\dd z\in\7P$. Clearly,
$\7P_{k}$ is the subspace of all $({k+1})$-homogeneous vector fields
in $\7P$ if $k\ge1$ and is $0$ otherwise.

Now assume that $F\subset F$ is a real analytic submanifold,
$M=F\oplus iV$ is the corresponding tube manifold and $a\in F$ is a
given point. For $\7g:=\hol(M,a)$ then put
$$\7g_{k}:=\7g\cap \7P_{k}\Steil{for all}k\in\ZZ\,.$$ Then $\{iv\dd
z:v\in V\}\subset\7g_{-1}$, and equality holds if $M$ is
holomorphically nondegenerate. Every $g\in\Aut(M,a)$ induces a Lie
algebra automorphism $\Theta:=g_{*}$ of $\7g$ in a canonical way: In
terms of local holomorphic representatives in suitable open
neighbourhoods of $a$ in $E$
$$\Theta(f(z)\dd z)=g'(g^{-1}z)(f(g^{-1}z))\Leqno{GP}$$ holds. This
implies immediately

\Lemma{GO} $g\mapsto g_{*}$ defines a group monomorphism
$\Aut(M,a)\;\hookrightarrow\;\Aut(\7g)$.  

\Proof Suppose that $\Theta:=g_{*}=\id$ for a $g\in\Aut(M,a)$. Then
$\Theta$ extends to a complex Lie automorphism of
$\7l:=\7g+i\7g\subset\hol(E,a)$ and leaves $\7P_{-1}\subset\7l$
element wise fixed. This implies that $g$ is represented by a
translation of $E$.  But $g(a)=a$ then gives $g=\id$.\qed

\KAP{cones}{Tube manifolds over cones}

In this section we always assume that the submanifold $F\subset V$ is
conical (that is, $x\in T_xF$ for every $x\in F$) and that $a\in F$ is
a given point. Then, for $M:=F\oplus iV$, the Lie algebra
$\7g:=\hol(M,a)$ contains the Euler vector field $\delta:=z\dd z$.

\Lemma{DC} Suppose that $\7g$ has finite
dimension. Then $\7g\subset\7P$ and for $\7g_{k}:=\7g\cap\7P_{k}$
$$\7g=\bigoplus\limits_{k\ge-1}\7g_{k}\,,\quad[\7g_{k},\7g_{l}]\subset
\7g_{k+l}\Steil{and}\7g_{-1}=\{iv\dd z: v\in V\}\,.\Leqno{AK}$$ In
particular, every $f(z)\dd z\in\7g$ is a polynomial vector field on
$E=V\oplus iV$ with $f(iV)\subset iV$.

\Proof Consider $\7l:=\7g\oplus\, i\7g\subset\hol(E,a)$, which
contains the vector field $\eta:=(z-a)\dd z$.  We first show
$\7l\subset \7P$: Fix an arbitrary $\xi:=f(z) \dd z\in\7l$. Then in a
certain neighbourhood of $a\in E$ there exists a unique expansion
$\xi=\sum_{k\in\NN}\xi_{k}$, where $\xi_{k}=p_{k}(z-a)\dd z$ for a
$k$-homogeneous polynomial map $p_{k}:E\to E$. It is easily verified
that the vector field $\ad(\eta)\xi\in\7l$ has the expansion
$\ad(\eta)\xi=\sum_{k\in\NN}(k{-}1)\xi_{k}$. Now assume that for
$d:=\dim(\7l)$ there exist indices $k_{0}<k_{1}<\dots<k_{d}$ such that
$\xi_{k_{l}}\ne0$ for $0\le l\le d$. Since the Vandermonde matrix
$\big((k_{l}-1)^{j}\big)$ in non-singular, we get that the vector
fields $(\ad\eta)^{j}\xi=\sum_{k\in\NN}(k{-}1)^{j}\xi_{k}$, $0\le j\le
d$, are linearly independent in $\7l$, a contradiction. This implies
$\xi\in\7P$ as claimed.\nline Since $\7g\subset\7P$ has finite
dimension, every $\xi\in\7g$ is a finite sum
$\xi=\sum_{k=-1}^{m}\xi_{k}$ with $\xi_{k}\in\7P_k$ and $m\in\NN$ not
depending on $\xi$. For every polynomial $p\in\RR[X]$ then
$p(\ad\delta)\xi1=\sum_{k=-1}^{m}p(k)\xi_{k}$ shows
$\xi_{k}\in\7g_{k}$ for all $k$, that is, $\7g=\oplus\7g_{k}$. The
identity $\7g_{-1}=\{iv\dd z: v\in V\}$ follows from the fact that
$\7g_{-1}$ is totally real in $\7P_{-1}$and this implies $f(iV)\subset
iV$ for all $f(z)\dd z\in\7g_{k}$ by
$[\7g_{-1},\7g_{k}]\subset\7g_{k-1}$ and induction on $k$.\qed

\Proposition{ZB} Assume that $\7g:=\hol(M,a)$ has finite dimension and
that $F'\subset V$ is a further conical submanifold with tube manifold
$M'=F'\oplus iV$. Then for every $a'\in F'$ every CR-isomorphism
$(M,a)\to(M',a')$ of manifold germs is rational.

\Proof Let $\7l:=\7g+\,i\7g$ and $\7l':=\7g'+\,i\7g'$ for
$\7g':=\hol(M',a')$. Fix a CR-isomorphism $(M,a)\to(M',a')$. This is
represented by a biholomorphic mapping $U\to U'$ with $g(a)=a'$ and
$g(U\cap M)=U'\cap M'$ for suitable connected open neighbourhoods
$U,U'$ of $a,a'\in E$. Then $g$ induces a Lie algebra isomorphism
$\Theta=g_{*}:\7l\to\7l'$, whose inverse is given by
$$\Theta^{-1}\big(f(z)\dd
z\re1\big)=g'(z)^{-1}\!f\big(g(z)\big)\re2\dd z\;.\Leqno{ZA}$$ Lemma
\ruf{AK} shows that $\7l$ and $\7l'$ consist of polynomial vector
fields. Consequently there exist {\sl polynomial} maps $p:E\to E$ and
$ q:E\to\End(E)$ such that
$$\Theta^{-1}\big(z\dd z\big)=p(z)\dd z\Steil{and}\Theta^{-1}\big(e\dd
z\big)=\big(q(z)e\big)\dd z$$ for all $e\in E$. Then \Ruf{ZA} implies
$g'(z)^{-1}=q(z)$ and $g'(z)^{-1}g(z)=p(z)$, that is,
$$g(z)=q(z)^{-1}p(z)$$ in a neighbourhood of $a\in E$.\qed

\Remark{} The tube $\5M$ over the future light cone shows that in
Proposition \ruf{LO} `rational' cannot always be replaced by `affine'.

\medskip Recall that $\aut(M,a)\subset\7g=\hol(M,a)$ is defined as the
isotropy subalgebra at $a$ and $\Aut(M,a)$ is the CR-automorphism group
of the manifold germ $(M,a)$, the so called {\sl stability group}.

\Proposition{LO} The following conditions are equivalent in case
$\7g=\hol(M,a)$ has finite dimension.  \0 $\7g_{1}=0$.  \1$\7g_{k}=0$
for all $k\ge1$.  \1 Every germ in $\Aut(M,a)$ can be represented by a
linear transformation $g\in\GL(V)\subset\GL(E)$ with $g(a)=a$.
\1 The tangential representation $\;h\mapsto h'(a)$ induces a group
monomorphism $\Aut(M,a)\;\hookrightarrow\;\GL(V)$.
\smallskip\noindent Each of these conditions is satisfied if
$\aut(M,a)=0$. On the other hand, if (iv) is satisfied then the image
of the group monomorphism is contained in the subgroup
$$\{g\in\GL(V):g\li1\7g_{0}=\7g_{0}\li1g,\;g(a)=
a\steil{and}g(T_{a}F)=T_{a}F \}\,,$$ where $\7g_{0}$ is considered in
the natural way as linear subspace of $\,\End(V)$.

\Proof Let $\7l:=\7g+\,i\7g\subset\hol(E,a)$ and
$\7l_{k}:=\7g_{k}+\,i\7g_{k}$ for all $k$.\nline \To12 A direct check
shows that, given any $\xi\in\7l$, the equality $[\7l_{-1},\xi]=0$
implies $\xi\in \7l_{-1}$.  Suppose $\7l_{k}\ne0$ for some $k>1$ and
that $k$ is minimal with respect to this property. Consequently,
$[\7l_{-1},\7l_{k}]\subset\7l_{k-1}=0$ would imply
$\7l_{k}\subset\7l_{-1}$, i.e., since $\7l_{-1}\cap\7l_{k}=0$, we
would have $\7l_{k}=0$. This contradicts our assumption $\7l_{k}\ne0$.
\nline\To23 Then $\eta:=(z-a)\dd z\in\7l$ and $\7l=\7l_{-1}\oplus\7k$,
where $\7k:=\{\xi\in\7l:\xi_{a}=0\}$ is the kernel of $\ad(\eta)$ in
$\7l$. Every germ in $\Aut(M,a)$ can be represented by a locally
biholomorphic map $ h:U\to E$ with $h(a)=a$, where $U$ is an open
neighbourhood of $a$ in $E$. The Lie algebra automorphism
$\Theta\in\Aut(\7l)$ induced by $h$ leaves the isotropy subalgebra
$\7k\subset\7l$ invariant. But $\eta$ is the unique vector field
$\xi\in\7k$ such that $\ad(\xi)$ induces the negative identity on the
factor space $\7l\!/\7k$. Therefore $\Theta(\eta)=\eta$ and hence also
$\Theta(\7l_{-1})=\7l_{-1}$ since $\7l_{-1}$ is the $(-1)$-eigenspace of
$\ad(\eta)$. From $\Theta(\7g_{-1})=\7g_{-1}$ we derive $h(z)=g(z)+c$
with $g\in\GL(V)$ and $c=a-g(a)\in V$. Taking the commutator of
$t\cd\id$ with $h$ we see that the translation $z\mapsto z+(t-1)c\,$
leaves $M$ invariant for all $t>0$. Differentiation by $t$ gives $c\dd
z\in\7g_{-1}$, that is, $c=0$ and $h=g$. \nline\To34 This is trivial.
\nline\To41 Let $\xi\in\7g_{1}$ be an arbitrary vector field. Then
there exists a unique symmetric bilinear map $b:E\times E\to E$ with
$\xi=b(z,z)\dd z$. Now $(\ad ia\dd z)^{2}\xi=-2b(a,a)\dd z\in\7g$,
that is, $\eta:=h(z)\dd z$ is in $\aut(M,a)$, where
$h(z):=b(z,z)-b(a,a)$. For every $t\in\RR$ therefore the
transformation $\psi_{t}:=\exp(t\eta)\in\Aut(M,a)$ has derivative
$\psi_{t}'(a)=\exp(th'(a))\in\GL(E)$ in $a$. But
$\psi_{t}'(a)\in\GL(V)$ by (iv) and thus $2b(a,v)=h'(a)v\in V$ for all
$v\in V$. On the other hand $b(a,v)\in iV$ by Lemma \ruf{ZI}, implying
$\psi_{t}'(a)=\id$ for all $t\in\RR$. By the injectivity of the
tangential representation therefore $\eta_{t}$ does not depend on $t$
and we get $\xi=0$. This proves the equivalence of (i) -- (iv).

Now suppose $\aut(M,a)=0$ and that there exists a non-zero
vector field $\xi\in\7g_{1}$. Then $\xi_{a}\in iV$ and there exists
$\eta\in\7g_{-1}$ with $\xi-\eta\in\aut(M,a)$, a contradiction.
Finally, assume (iv) and let $g\in\Aut(M,a)$ be arbitrary. By (iii)
then $g\in\GL(V)$ and $g(a)=a$. The induced Lie algebra automorphism
$\Theta$ maps every $\lambda\in\7g_{0}$ to $g\lambda g^{-1}\in \7g_0$,
compare \Ruf{ZA}.\qed

\bigskip We close the section with an example that illustrates how in
certain cases $\7g=\hol(M,a)$ can be explicitly computed with the
above results. We note that we do not know a single example with
$\dim\7g<\infty$ and $\dim\7g_{k}>\dim\7g_{-k}$ for some $k\in\NN$.

\medskip\Example{EB} Fix integers $p\ge q\ge0$ with $n:=p+q\ge3$ and a
real number $\alpha$ with $\alpha^2\ne\alpha$. Then for
${\RR_{+}:=\{t\in\RR:t>0\}}$
$$F=F_{p,q}^{\alpha}:=\Big\{x\in\RR_{+}^{n}:\sum_{j=1}^{n}\epsilon_{j}
x_{j}^{\alpha}= 0\Big\}\,,\quad\epsilon_{j}:=\cases{1&$j\le
p$\cr-1&otherwise\cr}\,,$$ is a hypersurface in
$V:=\RR^{n}$. Furthermore, $F$ is a cone and therfore $\dim
K_{a}F\ge1$ for every $a\in F$. On the other hand, the second
derivative at $a$ of the defining equation for $F$ gives a
non-degenerate symmetric bilinear form on $V\times V$, whose
restriction to $T_{a}F\times T_{a}F$ then has a kernel of dimension
$\le1$. Therefore $\dim(K_{a}F)=1$ for every $a\in F$ and by Corollary
\ruf{AF} the CR-manifold $M=M_{p,q}^{\alpha}:=F_{p,q}^{\alpha}\oplus
iV$ is everywhere $2$-nondegenerate (compare Example 4.2.5 in
\Lit{EBFT} for the special case $n=\alpha=3$). Since $M$ as
hypersurface is also minimal, $\7g=\hol(M,a)$ has finite dimension.

For the special case $\alpha=2$ and $q=1$ the above cone
$F=F_{n-1,1}^{2}$ is an open piece of the future light cone
$$\{x\in\RR^{n}:x_{1}^{2}+\dots+x_{n-1}^{2}=x_{n}^{2},\;x_{n}>0\}$$ in
$n$-dimensional space-time, which is affinely homogeneous. In
\Lit{KAZT} it has been shown that for the corresponding tube manifold
$M$ the Lie algebra $\7g=\hol(M,a)$ is isomorphic to $\7{so}(n,2)$ for
every $a\in M$. In case $q>1$ the following result is new:

\noindent {\bf Case $\Alpha\bf=2$:} Consider on $\CC^{n}$ the
symmetric bilinear form $\langle
z|w\rangle:=\sum\epsilon_{j}z_{j}w_{j}$. Then $F$ is an open piece of
the hypersurface $\{x\in\RR^{n}:x\ne0\,,\;\langle
x|x\rangle=0\}$, on which the reductive group $\RR^{*}\cd\O(p,q)$ acts
transitively. Therefore, $\7s_{0}^{}:=\RR\delta\oplus\7{so}(p,q)$ is
contained in $\7g_{0}$. One checks that
$$\7s:=\7g_{-1}\oplus\7s_{0}^{}\oplus\7s_{1}^{}\,,\qquad\7s_{1}^{}:=
\big\{\big(2i\langle c|z\rangle z-i\langle z|z\rangle c\big)\dd
z:c\in\RR^{n}\big\}$$ is a Lie subalgebra of $\7g$. The radical $\7r$
of $\7s$ is $\ad(\delta)$-invariant and hence of the form
$\7r=\7r_{-1}\oplus\7r_{0}\oplus\7r_{1}$ for
$\7r_{k}^{}:=\7r\cap\7g_{k}$. From $\7{so}(p,q)$ semisimple we
conclude $\7r_{0}\subset\RR\delta$. But $\delta$ cannot be in $\7r$
since otherwise $\7g_{-1}\subset\7r$ would give the false statement
$[\7g_{-1},\7s_{1}^{}]\subset\RR\delta$. Therefore $\7r_{0}^{}=0$, and
$[\7g_{-1},\7r_{1}^{}]=[\7r_{-1}^{},\7s_{1}^{}]=0$ implies
$\7r=0$. Now Proposition 3.8 in \Lit{KAZT} implies $\7g=\7s$, and, in
particular, that $\7g$ has dimension ${n+2\choose2}$. In fact, it can
be seen that $\7g$ is isomorphic to $\7{so}(p+1,q+1)$.

\medskip\noindent {\bf Case $\Alpha$ an integer $\bf\ge3$:} Then $F$ is
an open piece of the real-analytic submanifold
$$S:=\Big\{x\in\RR^{n}:x\ne0\Steil{and}\sum_{j=1}^{n}\epsilon_{j}
x_{j}^{\alpha}= 0\Big\}\Leqno{NP}$$ which is connected in case
$q>1$ and has two connected components otherwise. For every
$x\in\RR^{n}$ let $d(x)\in\NN$ be the cardinality of the set
$\{j:x_{j}=0\}$. It is easily seen that $\dim K_{x}S=1+d(x)$ holds for
every $x\in S$. Now consider the group
$$\GL(F):=\{g\in\GL(V):g(F)=F\}\,.$$ Every $g\in\GL(F)$ leaves $S$ and
hence also $H:=\{x\in\RR^{n}:d(x)>0\}$ invariant, that is, $g$ is the
product of a diagonal with a permutation matrix. Inspecting the action
of $\GL(F)$ on $\{c\in\overline F:d(c)=n-2\}$ we see that $\GL(F)$ as
group is generated by $\RR_{+}\!\cd\id$ and certain coordinate
permutations. As a consequence, $\7g_{0}=\RR\delta$ and, by a
simplified version of the argument following \Ruf{UP}, $\7g_{k}=0$ for
all $k>0$ is derived. In particular, $\dim\7g=n+1<\dim M$ for the tube
manifold $M=F\oplus iV$, that is, $M$ is not locally homogeneous. For
$n=3$ this gives an alternative proof for Proposition 6.36 in
\Lit{EBEN}.

\Proposition{} Let $a,a'\in F=F^{\alpha}_{p,q}$ be arbitrary
points. Then in case $3\le\alpha\in\NN$ the CR-manifold germs $(M,a)$
and $(M,a')$ are CR-equivalent if and only if $a'\in\GL(F)a$.

\Proof Suppose that $g:(M,a)\to(M,a')$ is an isomorphism of
CR-manifold germs. From $\7g'=\7g_{-1}'\oplus\,\RR\delta$ for
$\7g':=\hol(M',a')$ we conclude that $g$ is represented by a linear
transformation in $\GL(V)$ that we also denote by $g$. But then
$g(F)\subset S$ with $S$ defined in \Ruf{NP}. Because $g(F)$ has empty
intersection with $H$ we actually have $g(F)\subset F$. Replacing $g$
by its inverse we get the opposite inclusion, that is
$g\in\GL(F)$.\qed

\KAP{Some}{Levi degenerate CR-manifolds associated with an endomorphism}

Obviously, $2$ is the lowest CR-dimension for which there exist
homogeneous CR-manifolds that are Levi degenerate but not
holomorphically degenerate. In the following we present a large class
of such manifolds that contains all linearly
homogeneous conical tube manifolds of CR-dimension 2.

\Construction{CO}
Throughout this section let $V$ be a real vector space of dimension
$n\ge3$ and $d\ge1$ an integer with $d\le n-2$. Let furthermore
$\phi\in\End(V)$ be a fixed endomorphism and $\7h\subset\End(V)$ the
linear span of all powers $\phi^{k}$ for $0\le k\le d$. Then
$H:=\exp(\7h)\subset\GL(V)$ is an abelian subgroup and for given $a\in
V$ the orbit $F:=H(a)$ is a cone and an immersed submanifold of $V$
(not necessarily locally closed in case $n\ge4$). The tube
$M:=F\oplus iV\subset E$ is an immersed submanifold of $E$ and a
homogeneous CR-manifold in a natural way.\Formend

\Proposition{PR} Suppose that $M$ is minimal. Then $M$ is
$2$-nondegenerate and $a$ is a cyclic vector of $\phi$ (that is, $V$
is spanned by the vectors $\phi^{k}(a)$, $k\ge0$). Furthermore, $M$
has CR-dimension $d{+}1$, and the Levi kernel $K_{a}M$ has dimension
$1$.

\Proof Let $W\subset V$ be the linear span of all vectors
$\phi^{k}(a)$, $k\ge0$. Then $H\subset\RR[\phi]$ implies $H(a)\subset
W$ and hence $W=V$ by the minimality of $M$. Therefore, $a$ is a
cyclic vector, and the $\phi^{k}(a)$, $0\le k\le d$, form a basis of
the tangent space $T_{a}F$, that is, $M$ has CR-dimension $d+1$. Lemma
\ruf{CR} gives $\RR a\subset K_{a}F$ since $F$ is a cone in $V$. For
the proof of the opposite inclusion fix an arbitrary $w\in T_{a}F$
with $w\notin\RR a$. Then $w=\sum_0^m c_j \phi^j(a)$ with $c_m\ne0$
for some $m\ge1$, and $\phi^{d-m+1}(w)\notin T_{a}F$ shows $w\notin
K_{a}F$ by Proposition \ruf{US}. Therefore, $M$ is $2$-nondegenerate
by Corollary \ruf{AF}.\qed

From linear algebra it is clear that $\phi$ has a cyclic
vector if and only if minimal and characteristic polynomial of $\phi$
coincide. Furthermore, if $a,b\in V$ both are cyclic vectors of $\phi$
then there exists a transformation $g\in\RR[\phi]\subset\End(V)$ with
$b=g(a)$. But $g$ commutes with every element of the group
$H=\exp(\7h)$ and hence maps the orbit $H(a)$ onto $H(b)$. Since
$H(b)$ is not contained in a hyperplane of $V$, necessarily $g$ is in
$\GL(V)$, i.e. all cyclic vectors of $\phi$ give locally CR-equivalent
tube manifolds.

For $M$ in Proposition \ruf{PR} the Lie algebra $\7g_{0}$ contains the
$(d{+}1)$-dimensional Lie algebra $\7h$. We claim, that actually for
$\phi$ in `general position' $\7g_{0}=\7h$ holds. Indeed, suppose that
$\phi$ is diagonalizable with eigenvalues
$\alpha_{1},\dots,\alpha_{n}\in\CC$. For $N:=\{1,2,.\dots,n\}$ and
every subset $I\subset N$ of cardinality $\#I=d{+}1$ put
$$\Delta_{I}:=\big\{(\alpha_{k}^{}-\alpha_{j}^{},\alpha_{k}^{2}-
\alpha_{j}^{2}, \dots,\alpha_{k}^{d}-\alpha_{j}^{d})\in\CC^{d}:k\in
N,\,j\in I\big\}\,,$$ which obviously contains the origin
$0\in\CC^{d}$.

\Proposition{UH} Suppose that $M$ is minimal, $\phi$ is diagonalizable
and that
$$\bigcap\limits_{\#I=d+1}\!\Delta_{I}=\{0\}\;.\Leqno{DU}$$ Then
$\7g=\7g_{-1}\oplus\7h$ and $\aut(M,a)=0$.

\Proof By the minimality of $M$ the eigenvalues
$\alpha_{1},\alpha_{2},\dots,\alpha_{n}$ are mutually distinct. The
Lie algebra $\7g_{0}$ consists of all operators $\mu\in\End(V)$ that
are tangent to $F=H(a)$, or equivalently, such that for every
$t=(t_{1},t_{2},\dots,t_{d})\in\RR^{d}$ there exist real coefficients
$r_{0},r_{1},\dots,r_{d}$ with
$$\mu(e^{t_{1}\phi+t_{2}\phi^{2}+\dots+t_{d}\phi^{d}}a)=
(r_{0}+r_{1}\phi+\dots+r_{d}\phi^{d})e^{t_{1}\phi+t_{2}\phi^{2}
+\dots+t_{d}\phi^{d}}a \,.\Leqno{PS}$$ Since the orbit $H(a)$ contains
a basis of $V$ the matrix $\mu$ is uniquely determined by the function
tuple $r_{0},\dots,r_{d}$ on $\RR^{d}$. Going to the complexification
$E$ of $V$ we may assume that $E=\CC^{n}$ and that $\phi$ is given by
the diagonal matrix with diagonal entries $\alpha_1,\dots,\alpha_n$.
Since $a=(a_{1},\dots,a_{n})\in\CC^{n}$ is a cyclic vector of $\phi$
we have $a_{j}\ne0$ for every $j$. Then \Ruf{PS} for
$\mu=(\mu_{jk})\in\CC^{n\Times n}$ reads
$$r_{0}^{}+\alpha_{j}r_{1}^{}+\alpha_{j}^{2}r^{}_{2}+\dots+
\alpha_{j}^{d}r^{}_{d} =a_{j}^{-1}\sum_{k=1}^{n}\mu_{jk}^{}a_{k}^{}e^{
(\alpha_{k} -\alpha_{j})t_{1} +(\alpha_{k}^{2}
-\alpha_{j}^{2})t_{2}+\dots+(\alpha_{k}^{d}-
\alpha_{j}^{d})t_{d}}\Leqno{EQ}$$ for all $j\in N$. For every subset
$I\subset N$ with $\#I=d{+}1$ the system of all linear equations
\Ruf{EQ} with $j\in I$ has a unique solution for $r_{0},\dots,r_{d}$
since the coefficient matrix is of Vandermonde type. In particular,
every $r_{p}$ is a complex linear combination of functions
$e^{\beta_{1}t_{1}+\beta_{2}t_{2}+\dots+\beta_{d}t_{d}}$ with
$(\beta_1,\dots,\beta_{d})\in\Delta_{I}$. Note that given any set of
pairwise different vectors $(\beta_1,\ldots,\beta_d)$ the
corresponding functions
$e^{\beta_{1}t_{1}+\beta_{2}t_{2}+\dots+\beta_{d}t_{d}}$ are linearly
independent. Consequently, if follows from \Ruf{DU} that actually
every $r_{p}$ is constant on $\RR^{d}$. But this implies
$\dim\7g_{0}\le d+1=\dim\7h$ and thus $\7g_{0}=\7h$. \nline For the
second claim let $\xi\in\7g_{1}$ be an arbitrary vector field. With
respect to $E=\CC^{n}$ and $\phi$ the diagonal matrix as before, we
may write
$$\xi=\sum_{j,k,p=1}^{n}c_{p}^{jk}z_{j}z_{k}\dd{z_{p}}\,,\qquad
c_{p}^{jk}=c_{p}^{kj}\in\CC\,.\Leqno{UP}$$ Then
$$\xi_{r}:=\big[\dd{z_{r}},\xi\big]=2\sum_{k,p}c^{kr}_{p}z_{k}
\dd{z_{p}}\in \7g_0\oplus i\7g_0=\7h\oplus i\7h$$ for all $r\in N$
implies $c_{p}^{kr}=0$ if $k\ne p$ and thus, by the symmetry in $k$
and $r$, we have $\xi_{r}=2c^{rr}_{r}z_{r}\dd{z_{r}}$. On the other
hand, if $\lambda(z)\dd z\in \7h\oplus i\7h$ is non-zero then
$\rank(\lambda)\ge n-d\ge 2$. But this implies $\xi_{r}=0$ for all
$r$, and in turn $\xi=0$, i.e. $\7g_{1}=0$. By Proposition \ruf{LO}
therefore $\7g_{k}=0$ holds for all $k>0$. Finally,
$\dim(\7g)=n+d+1=\dim(M)$ implies $\aut(M,a)=0$.\qed

\smallskip Notice that in case $n=\dim(V)=3$ the condition \Ruf{DU} is
equivalent to
$$\alpha_{1}\ne(\alpha_{2}+\alpha_{3})/2,\quad\alpha_{2}\ne
(\alpha_{3}+\alpha_{1})/2\Steil{and}\alpha_{3}\ne
(\alpha_{1}+\alpha_{2})/2\,.$$ This will be applied in examples
\ruf{EY} and \ruf{EX} (where \ruf{DU} holds) and in example \ruf{EI}
(where \ruf{DU} fails).

\medskip The conclusion in Proposition \ruf{UH} also holds in many
cases where $\phi$ is not diagonalizable, see e.g. Example \ruf{EZ}
below. In any case, we get as a consequence of the above

\Proposition{UQ} Suppose that the tube $M$ over the cone $F=H(a)$ is
simply connected and that $\Aut(M,a)$ is the trivial group. Then the
following properties hold: \0 Let $M'$ be an arbitrary homogeneous
CR-manifold and $D\subset M$, $D'\subset M'$ non-empty domains. Then
every real-analytic CR-isomorphism $h:D\to D'$ extends to a
real-analytic CR-\-isomorphism $M\to M'$. \1 Let $M'$ be an arbitrary
locally homogeneous CR-manifold and $D'\subset M'$ a domain that is
CR-equivalent to $M$. Then $D'=M'$.

\Proof The group $G:=\Aut(M)$ acts simply transitive on $M$ and has
trivial center by Lemma \ruf{ZI}. \nline {\bf ad (i)} We may assume
$a\in D$. To every $g\in G$ with $g(a)\in D$ there exists a
transformation $g'\in G':=\Aut(M')$ with $hg(a)=g'h(a)$. Because of
$\Aut(M',a')=\{\id\}$ the transformation $g'$ is uniquely determined
by $g$ and satisfies $hg=g'h$ in a neighbourhood of $a$. Since the Lie
group $G$ is simply connected $g\mapsto g'$ extends to a group
homomorphism $G\to G'$ and $h$ extends to a CR-covering map $h:M\to
M'$. The deck transformation group $\Gamma:=\{g\in G:gh=h\}$ is in the
center of $G$, which is trivial by Lemma \ruf{ZI}. Therefore, $h:M\to
M'$ is a CR-isomorphism. \nline {\bf ad (ii)} The proof is essentially
the same as for Proposition 6.3 in \Lit{FEKA}.\qed

The condition 'locally homogeneous' in Proposition \ruf{UQ}.ii cannot
be omitted. A counter example is given for every integer $\alpha\ge3$
by the tube $M'\subset\CC^{3}$ over
$$F':=\{x\in\RR^{3}:x_{3}=x_{1}(x_{2}/x_{1})^{\alpha},\;x_{1}>0\}\,.$$
Then the tube $M$ over $F:=\{x\in F':x_{2}>0\}$ is the Example
\ruf{EX} below that satisfies the assumption of Proposition \ruf{UQ}.

The manifold $M$ in Proposition \ruf{UQ} is simply connected, for
instance, if there exist eigenvalues
$\alpha_{0},\alpha_{1},\dots,\alpha_{d}$ of $\phi$ such that
$\det(A+\overline A)\ne0$, where $A=(\alpha_{j}^{k})_{0\le j,k\le d}$
is the corresponding Vandermonde matrix. A criterion for the
triviality of the stability group is given by the following statement.

\Proposition{EO} Suppose that $a\in V$ is a cyclic vector of $\phi\in
\End V$ and that $\7g_{1}=0$ for $\7g=\hol(M,a)$, $M=F\oplus iV$ and
$F=H(a)\subset V$ as above. Then $\Aut(M,a)=\{\id\}$ and
$\7g_{0}=\7h$ if
$$\{g\phi g^{-1}:g\in\GL(V)\steil{with}g(a)=a\}\;\cap\;\7g_{0}=
\{\phi\}\,.$$

\Proof Fix $g\in\Aut(M,a)$. By Proposition \ruf{LO}, $g\in\GL(V)$,
$g(a)=a$, $g\phi g^{-1}\in\7g_0$ and hence, by our assumption,
$g\phi=\phi g$. But then $g\phi^{k}(a)=\phi^{k}g(a)=\phi^{k}(a)$ for
every $k$ implies $g=\id$ since $a$ is a cyclic vector. For the second
claim note that by construction of $M$ the restriction of the
evaluation map $\epsilon_a:\7g\to T_aM\;\cong\;\7g_{-1}\oplus \7h$ is
surjective.  Since $\Aut(M,a)=\{\id\}$ implies $\aut(M,a)=0=\ker
(\epsilon_a)$, the claim follows.\qed

\KAP{Examples}{Homogeneous 2-nondegenerate CR-manifolds in dimension 5}

\medskip In the following we present (up to local linear equivalence)
all linearly homogeneous cones $F\subset\RR^{3}$, for which $M=F\oplus
i\RR^3$ is $2$-nondegenerate (that is, $F$ is not contained in a
hyperplane of $\RR^{3}$). All these are obtained by the recipe of the
preceding section: Choose a linear operator $\phi$ on $\RR^{3}$ having
a cyclic vector $a\in\RR^{3}$. Then $\7h:=\RR\id+\RR\phi$ is an
abelian Lie algebra and for the corresponding connected Lie subgroup
$H\hookrightarrow\GL(3,\RR)$ the orbit $F:=H(a)$ is a homogeneous cone
with $M:=F\oplus i\RR^3$ Levi degenerate. In case $\phi'$ is another
choice such that $\7h'=\RR\id\oplus\RR\phi'$ is conjugate in
$\GL(3,\RR)$ to $\7h$, then the spectra of $\phi$, $\phi'$ differ in
$\CC$ only by an affine transformation $z\mapsto rz+s$ with
$r,s\in\RR$ and $r\ne0$.

\Example{EI} Let
$F=F_{0}:=\{x\in\RR^{3}:x_{1}^{2}+x_{2}^{2}=x^{2}_{3},\,x_{3}>0\}$ be
the future light cone. This occurs for the choice
$\phi:=x_{2}\dd{x_{1}}-x_{1}\dd{x_{2}}$ having spectrum $\{\pm i,0\}$.
The vector field $\phi':=x_{3}\dd{x_{2}}+x_{2}\dd{x_{3}}$ has spectrum
$\{\pm 1,0\}$ and generates with $\delta$ a $2$-dimensional Lie
algebra $\7h'$ that is tangent to $F$. An $\7h'$-orbit is the domain
$F':=\{x\in F:x_{1}>0\}$ in $F$. The same happens with the nilpotent
vector field
$$\phi'':=\phi+\phi'=x_{2}\dd{x_{1}}+
(x_{3}-x_{1})\dd{x_{2}}+x_{2}\dd{x_{3}}\,,$$ that generates with
$\delta$ a $2$-dimensional Lie algebra $\7h''$. Then $F'':=\{x\in
F:x_{2}<x_{3}\}$ is the unique open $\7h''$-orbit in $F$. For the tube
$M=F\oplus i\RR^{3}$ it is known that $\7g:=\hol(M,a)$ is a
$10$-dimensional simple Lie algebra, isomorphic to $\7{so}(2,3)$,
compare \Lit{KAZT}, \Lit{FEKA}. The convex hull $\hat\, F$ of $F$ is
$\hat\, F=\{x\in\RR^{3}:x^{2}_{3}\ge x_{1}^{2}+x_{2}^{2},\,x_{3}>0\}$.

\Example{EY} For $\alpha>0$ let $F=F_{\alpha}\subset\RR^{3}$ be the
orbit of $(1,0,1)$ under the group of all linear transformations
$x\mapsto e^{s}(\cos t\,x_{1}-\sin t\,x_{2},\sin t\,x_{1}+\cos
t\,x_{2},e^{\alpha t}x_{3})$, $s,t\in\RR$. With
$r:=(x_{1}^{2}+x_{2}^{2})^{1/2}$, the manifold $F$ is given in
$\{x\in\RR^{3}:r>0\}$ by the explicit equations
$$x_{3}=r\exp\big(\alpha\cos^{-1}(x_{1}/r)\big)=r\exp\big
(\alpha\sin^{-1}(x_{2}/r)\big)\,,\leqno{(*)}$$ where locally always
one of these suffices. A suitable choice is
$\phi=x_{1}\dd{x_{2}}-x_{2}\dd{x_{1}}+\alpha x_{3}\dd{x_{3}}$ with
spectrum $\{\pm i,\alpha\}$. The convex hull of $F$ is $\hat\,
F=\{x\in\RR^{3}:x_{3}>0\}$. Notice that replacing $\alpha$ by
$-\alpha$ would lead to a linearly equivalent cone.

\medskip   

\Example{EZ} Let $F=F_{\infty}\subset\RR^{3}$ be the
orbit of $(1,0,1)$ under the group of all linear transformations
$x\mapsto e^{s}(x_{1},x_{2}+tx_{1},e^{t}x_{3})$ with
$s,t\in\RR$, that is,
$$F=\big\{x\in\RR^{3}:x_{1}>0,\;x_{3}^{}= x_{1}^{}\re1e^{ x_{2}/x_{1}}
\big\}\,.$$ Here $\phi=x_{1}\dd{x_{2}}+x_{3}\dd{x_{3}}$ has
spectrum $0,0,1$ and is not diagonalizable.

\medskip   

\Example{EX} For $\alpha<-1$ let
$$F=F_{\alpha}:=\{x\in\RR^{3}_{+}:x_{3}=x_{1}(x_{2}/x_{1})^{\alpha}\}\,.$$
Here $\phi=x_{2}\dd{x_{2}}+\alpha x_{3}\dd{x_{3}}$ has distinct real
eigenvalues $\{\alpha,0,1\}$. Notice that the limit case $\alpha=-1$
has already been discussed in \ruf{EI}.

\bigskip

For every $F$ and $M=F\oplus i\RR^3$ in the examples \ruf{EY} --
\ruf{EX} consider for $\7g:=\hol(M,a)$ the decomposition \Ruf{AK}. By
a straightforward calculation it is seen that always
$\7g_{0}=\RR\delta\oplus\RR\phi$ and $\7g_{k}=0$ for $k>0$ (for
\ruf{EY} and \ruf{EX} this also follows from Proposition \ruf{UH} and
for \ruf{EZ} use Proposition \ruf{EO}). In particular, $\7g$ is a
solvable Lie algebra of dimension $5$ with commutator
$[\7g,\7g]=\7g_{-1}$ of dimension $3$. Including now also Example
\ruf{EI} let us for better distinction write $F=F_{\alpha}$,
$M=M_{\alpha}$ and $a=a_{\alpha}$, where $\alpha\in\RR\cup\{\infty\}$
satisfies $\alpha<-1$ or $\;\alpha\ge0$.

\Lemma{LM} The real Lie algebras $\hol(M_{\alpha},a_{\alpha})$ are
mutually non-isomorphic.

\Proof Since $\5M=M_{0}$ is the only manifold among all $M_{\alpha}$
such that $\hol(M_{\alpha},a_{\alpha})$ is not of dimension $5$ we
restrict our attention to $M=M_{\alpha}$ and $a=a_{\alpha}$ with
$\alpha\ne0$.  Then, as noted before, $\7g_{-1}=[\7g,\7g]$ for
$\7g=\hol(M,a)$, and $\7g_{-1}\oplus\;\RR\delta$ is the set of all
$\xi\in\7g$ such that the restriction
$\rho(\xi):=\ad(\xi)|_{\7g_{-1}}$ acts as a multiple of the identity
on $\7g_{-1}$. Therefore its complement
$\7g\!\backslash(\7g_{-1}\oplus\;\RR\delta)$ is an invariant of
$\7g$. \nline Denote by $\Gamma$ the group of all real affine
transformations $z\mapsto rz+t$ with $r,t\in\RR$ and $r\ne0$. Then
$\Gamma$ also acts in a natural way on the space $\CC^{3}/\7S_{3}$ of
all unordered complex number triples. For every
$\xi\;\in\;\7g\!\backslash(\7g_{-1}\oplus\;\RR\delta)$ the operator
$\rho(\xi)$ has precisely 3 (not necessarily distinct) complex
eigenvalues and hence determines a point in $\CC^{3}/\7S_{3}$. It is
easy to see that the set $\Sigma_{\alpha}$ of all points in
$\CC^{3}/\7S_{3}$ obtained in this way by elements
$\xi\notin(\7g_{-1}\oplus\;\RR\delta)$ is a $\Gamma$-orbit and also is
an invariant of the Lie algebra
$\7g=\hol(M_{\alpha},a_{\alpha})$. Since the map
$\alpha\mapsto\Sigma_{\alpha}$ is injective, the claim follows.\qed

\medskip The surfaces $F=F_{\alpha}\subset\RR^{3}$ in examples
\ruf{EI} -- \ruf{EX} are all cones. Cones belong to the class of
(locally) ruled surfaces and are characterized by the property that
all rules meet in a point. Another case of a ruled surface is if all
rules are the tangents of a fixed curve in $\RR^{3}$. By \Lit{EAEZ},
\p45, there is up to affine equivalence only one homogeneous ruled
surface of this type that can be described as follows.

\Example{EV} Let $S$ the union of all tangents to the {\sl twisted
cubic} $C:=\{(s,s^{2},s^{3})\in\RR^{3}:s\in\RR\}$, that is,
$$S=\{x\in\RR^{3}:x_{3}^2+4x_{2}^{3}-6x_{1}x_{2}x_{3}-
3x_{1}^{2}x_{2}^{2}+4x_{1}^{3}x_{3}=0\}\,.$$ The two one-parameter
groups
$$x\mapsto(e^{t}x_{1},e^{2t}x_{2},e^{3t}x_{3})\,,\qquad
x\mapsto(x_{1}+t,x_{2}+2tx_{1}+t^{2},x_{3}+3tx_{2}+3t^{2}x_{1}+t^{3})$$
leave $C$ and hence also $S$ invariant. They generate a connected
2-dimensional group $H$ of affine transformations that has three
orbits in $S$. Besides $C$ there are two further orbits in $S$ that
are interchanged by the involution $x\mapsto
(-x_{1},x_{2},-x_{3})$. Denote by $F$ the $H$-orbit containing
$a:=(1,0,0)$. Then $F$ is an affinely homogeneous surface in
$\RR^{3}$. Clearly, the Lie algebra $\7h$ of $H$ is generated by the
affine vector fields
$$\zeta=x_{1}\dd{x_{1}}+2x_{2}\dd{x_{2}}+3x_{3}\dd{x_{3}}\,,\qquad\quad
\eta=\dd{x_{1}}+2x_{1}\dd{x_{2}}+3x_{2}\dd{x_{3}}\,.$$ Here,
$T_{a}F=\{c\in\RR^{3}:c_{3}=0\}$ and
$K_{a}F=\{c\in\RR^{3}:c_{2}=c_{3}=0\}$ are easily checked with Lemma
\ruf{CR}. With Proposition \ruf{US} we see that $F$ is
$2$-nondegenerate. In particular, $\7g:=\hol(M,a)$ has finite
dimension.

For every integer $k$ denote by $\7g^ {(k)}$ the $k$-eigenspace of
$\ad(\zeta)$ in $\7g$. Then we have the $\ZZ$-gradation
$$\7g=\bigoplus_{k\ge-3}\7g^{(k)}\,.\Leqno{RU}$$ It can be seen that
$\7g^{(-3)}=\RR i\dd{z_{3}}\,$, $\;\7g^{(-2)}=\RR i\dd{z_{2}}\,$,
$\;\7g^{(-1)}=\RR i\dd{z_{1}}\oplus\RR\eta\,$,
$\;\7g^{(0)}=\RR\zeta\,$ and $\;\7g^{(k)}=0$ for $k>0$ (here $\zeta$
and $\eta$ are in the unique way extended to complex affine vector
fields on $\CC^{3}$). This implies that $\hol(M,a)$ is a solvable Lie
algebra of dimension $5$ with commutator ideal
$[\7g,\7g]=\bigoplus_{k<0}\7g^{(k)}$ of dimension $4$. On the other
hand, for $M=F\oplus i\RR^{3}$ with $F\subset\RR^{3}$ a cone from
Example \ruf{EI} -- \ruf{EX} the commutator of $\hol(M,a)$ has either
dimension $10$ (Example \ruf{EI}) or dimension $3$ (all the others). As
a consequence of Lemma \ruf{LM} we therefore get:

\Proposition{} The manifolds $M=F\oplus i\RR^{3}$, $F$ a cone from
examples \ruf{EI} -- \ruf{EV}, are all homogeneous $2$-nondegenerate
CR-manifolds. Furthermore, they are mutually locally
CR-inequivalent.\Formend

\Lemma{LS} Let $M:=F\oplus i\RR^{3}$ for $F=H(a)\subset\RR^{3}$ as in
Example \ruf{EV}. Then $\Aut(M,a)$ is the trivial group.

\Proof Let $\tilde M$ be the image of $M$ under the translation
$z\mapsto z-a$, where $a=(1,0,0)$. It is enough to show that
$\Aut(\tilde M,0)$ is the trivial group. Let
$\7l:=\tilde{\7g}+i\tilde{\7g}$ for $\tilde{\7g}:=\hol(\tilde M,0)$ and
put $\xi_{j}:=\dd{z_{j}}$ for $j=1,2,3$. From \Ruf{RU} and the
explicit description of all $\7g^{(k)}$ we see that
$\7l=\7P_{-1}\oplus\7l_{0}$ for $\7l_{0}=\CC\zeta\oplus\CC\eta$ with
$$\zeta:=z_{1}\dd{z_{1}}+2z_{2}\dd{z_{2}}+3z_{3}\dd{z_{3}},\quad
\eta:=2z_{1}\dd{z_{2}}+3z_{2}\dd{z_{3}}\,.$$ Now fix $g\in\Aut(\tilde
M,0)$ and denote by $\Theta=g_{*}$ the induced Lie algebra
automorphism of $\7l$.  Then $\7l_{0}$ is $\Theta$-invariant as
isotropy subalgebra at the origin. Further $\Theta$-invariant linear
subspaces are $\7n:=[\7l,\7l]$ and thus also
$$\CC\eta=\7l_{0}\cap\;\7n,\;\;\CC\xi_{1}=K_{0}\tilde
M,\;\;\CC\xi_{2}=H_{0}\tilde
M\cap[\7n,\7n],\;\;\CC\xi_{3}=[\CC\xi_{2},\CC\eta]\;.$$ As a
consequence, there exist $s_{1},s_{2},s_{3},u,v\in\CC^{*}$ and
$w\in\CC$ such that $\Theta(\xi_{j})=s_{j}\xi_{j},\Theta(\eta)=u\eta$
and $\Theta(\zeta)=v\zeta+w\eta$.  Then
$$s_{1}\xi_{1}=\Theta(\xi_{1})=\Theta([\xi_{1},\zeta])=[s_{1}\xi_{1},
v\zeta+w\eta]=s_{1}v\xi_{1}+2s_{1}w\xi_{2}$$ implies $v=1$, $w=0$ and
hence $\Theta(\zeta)=\zeta$. Since $s_{j}i\xi_{j}$ is contained in
$\tilde{\7g}$ we get $s_{j}\in\RR$ for $j=1,2,3$. Also the vector
fields $\zeta+\xi_{1}$ and $\eta+\xi_{1}+2\xi_{2}$ are contained in
$\tilde{\7g}$ and have $\Theta$-image in $\tilde{\7g}$, that is,
$(1-s_{1})\xi_{1}$ and $(1-u)\eta+2(1-s_{2})\xi_{2}$ are contained in
$\tilde{\7g}$.  This implies $s_{1}=s_{2}=u=1$. Finally, $\Theta$
applied to $[\xi_{2},\eta]=3\xi_{3}$ gives $s_{3}=1$. Therefore
$\Theta=\id$ and Lemma \ruf{GO} completes the proof.\qed

For the tube $\5M$ over the future light cone (that is Example
\ruf{EI}) there exist infinitely many homogeneous CR-manifolds that
are mutually globally CR-inequivalent but are all locally
CR-equivalent to $\5M$, compare \Lit{KAZT}. In contrast to this,
we have for the remaining examples:

\Proposition{} Let $M=F\oplus iV$ with $F$ from one of the examples
\ruf{EY} -- \ruf{EV}. Then $M$ is simply connected and $\Aut(M)$ is a
solvable Lie group of dimension 5 acting transitively and freely on
$M$. For every $a\in M$ the stability group $\Aut(M,a)$ is trivial and
every homogeneous real-analytic CR-manifold $M'$, that is locally
CR-equivalent to $M$, is already globally CR-equivalent to $M$.

\Proof It is easily checked that $F$ and hence $M$ is simply
connected. In case $F$ is a cone, the assumption in Proposition
\ruf{EO} is satisfied and the claim follows with Proposition \ruf{UQ}.
Therefore we may assume that $F$ is the submanifold of Example
\ruf{EV}. But then $\Aut(M,a)$ is the trivial group by Lemma \ruf{LS}
and $\Aut(M)$ has trivial center by Proposition \ruf{ZI}. But then
the proof of Proposition \ruf{UQ}
also applies to this case.\qed

\medskip With an argument from \Lit{FEKA} together with 2.5.10 in
\Lit{HORM} it can be seen that every continuous CR-function on
$M=F\oplus iV$, $F\subset V$ a cone as above, has a unique continuous
extension to the convex hull $\hat{M}=\hat{\re2F}\oplus i\RR^{3}$ of
$M$ which is holomorphic on the interior of $\hat{M}$ with respect to
$\CC^{3}$. In particular, if $F$ belongs to Example \ruf{EY}, then
every continuous CR-function on $M$ is real-analytic.

{\Klein 
\parindent 15pt\bigskip\bigskip{\noindent\gross References}
\bigskip 
\def\Springer{Ber\-lin-Hei\-del\-berg-New York: Sprin\-ger~}

\font\klCC=txmia scaled 900

\Ref{ANHI}{Andreotti, A., Hill, C.D.:} Complex characteristic coordinates and tangential Cauchy-Riemann equations. Ann. Scuola Norm. Sup. Pisa {\bf 26} (1972), 299-324.
\Ref{BERO}Baouendi, M.S., Ebenfelt, P., Rothschild, L.P.: {\sl Real Submanifolds in Complex Spaces and Their Mappings}. Princeton Math. Series {\bf 47}, Princeton Univ. Press, 1998.
\Ref{BAHR}{Baouendi, M.S., Huang, X., Rothschild, L.P.:} Regularity of CR mappings between algebraic hypersurfaces. Invent. Math. {\bf 125} (1996), 13-36.
\Ref{BUSH}Burns, D., Shnider, S.: Spherical hypersurfaces in complex manifolds. Invent. Math. {\bf 33} (1976), 223--246. 
\Ref{CHMO}Chern, S.S., Moser, J.K.: Real hypersurfaces in complex manifolds. Acta. Math. {\bf 133} (1974), 219-271.
\Ref{DOKR}Doubrov, B., Komrakov, B., Rabinovich, M.: Homogeneous surfaces in the three-dimenaional affine geometry. In: {\sl Geometry and Topology of Submanifolds, VIII,} World Scientific, Singapore 1996, pp. 168-178
\Ref{EAEZ}Eastwood, M., Ezhov, V.: On Affine Normal Forms and a Classification of Homogeneous Surfaces in Affine Three-Space. Geometria Dedicata {\bf 77} (1999), 11-69.
\Ref{EBFT}Ebenfelt, P.: Normal Forms and Biholomorphic Equivalence of Real Hypersurfaces in \hbox{\klCC\char131}$^3$. Indiana J. Math. {\bf 47} (1998), 311-366.
\Ref{EBEN}Ebenfelt, P.: Uniformly Levi degenerate CR manifolds: the 5-dimensional case. Duke Math. J. {\bf 110} (2001), 37--80.
\Ref{FEKA}Fels, G., Kaup, W.: CR-manifolds of dimension 5: A Lie algebra approach. J. Reine Angew. Math, to appear. {\ninett http://arxiv.org/pdf/math.DS/0508011}
\Ref{FEKP}Fels, G., Kaup, W.: Classification of locally homogeneous Levi degenerate CR-manifolds in dimension 5. In preparation.
\Ref{GAME}Gaussier, H., Merker, J.: A new example of a uniformly Levi degenerate hypersurface in \hbox{\klCC\char131}$^3$. Ark. Mat. {\bf 41} (2003), 85--94.
\Ref{HORM}H\"ormander,~L.: {\sl An Introduction to Complex Analysis in Several Variables}. Third edition. North-Holland Publishing Co. Amsterdam - New York, 1990.
\Ref{KAZT}Kaup, W., Zaitsev, D.: On local CR-transformations of Levi-degenerate group orbits in compact Hermitian symmetric spaces. J. Eur. Math. Soc.. To appear 
\Ref{SAKA} {Sakai, T.:} Riemannian geometry. Translations of Mathematical Monographs, 149. American Mathematical Society, Providence, RI, 1996.
\Ref{SEVL}Sergeev, A.G., Vladimirov, V.S.: Complex analysis in the future tube. In: {\sl Several Complex Variables II,} Encyclopedia Math. Sci {\bf 8}, \Springer 1994.
\Ref{TANA}Tanaka, N.: On the pseudo-conformal geometry of hypersurfaces of the space of $n$ complex variables. J. Math. Soc. Japan {\bf 14} (1962), 397-429.
\bigskip

\bigskip\centerline{Mathematisches Institut, Universit\"at
T\"ubingen, Auf der Morgenstelle 10, 72076 T\"ubingen, Germany}
\smallskip
\centerline{e-mail: {\ninett gfels@uni-tuebingen.de, kaup@uni-tuebingen.de}}
}

\closeout\aux\bye